\documentclass[11pt]{article}\usepackage[margin=1in]{geometry}

\usepackage[utf8]{inputenc}
\usepackage{amsmath,amssymb,amsthm,amsfonts}
\usepackage{mathtools}
\usepackage{enumitem}
\usepackage{booktabs}
\usepackage{array}
\usepackage{tikz}
\usetikzlibrary{automata,positioning,arrows.meta,calc}
\usepackage{pgfplots}
\pgfplotsset{compat=1.16}
\usepackage{hyperref}
\hypersetup{
    colorlinks=true,
    linkcolor=blue,
    filecolor=blue,
    citecolor=blue,
    urlcolor=blue
}
\usepackage{cleveref}
\usepackage{float}
\usepackage{caption}

\theoremstyle{plain}
\newtheorem{theorem}{Theorem}[section]
\newtheorem{lemma}[theorem]{Lemma}
\newtheorem{proposition}[theorem]{Proposition}
\newtheorem{corollary}[theorem]{Corollary}
\newtheorem{conjecture}[theorem]{Conjecture}
\theoremstyle{definition}

\newtheorem{remark}[theorem]{Remark}

\newcommand{\Z}{\mathbb{Z}}
\newcommand{\N}{\mathbb{N}}
\newcommand{\Cmap}{C}
\newcommand{\Prob}{\mathbb{P}}

\newcommand{\fa}{\mathfrak{a}}
\newcommand{\fb}{\mathfrak{b}}
\newcommand{\fc}{\mathfrak{c}}
\newcommand{\fd}{\mathfrak{d}}
\newcommand{\fe}{\mathfrak{e}}
\newcommand{\ff}{\mathfrak{f}}
\newcommand{\floor}[1]{\left\lfloor #1 \right\rfloor}
\newcommand{\TV}{\mathrm{TV}}
\DeclareMathOperator{\Freq}{Freq}
\DeclareMathOperator{\Bin}{Bin}
\DeclareMathOperator{\Pred}{Pred}

\title{Residue Automata Modulo $3^k$, Syracuse Random Variables, and Terminal-Branch Bias in Collatz Trajectories}

\author{
Sawon Pratiher$^{\ast}$
\\[6pt]
}

\date{}

\begin{document}

\maketitle

\begin{abstract}
We study the residues modulo $3^k$ visited by trajectories of the accelerated Collatz map $T(n)=n/2$ ($n$ even), $T(n)=(3n+1)/2$ ($n$ odd), starting from the six ``recurrent forms'' $9n+r$, $r\in\{1,2,4,5,7,8\}$. (i)~We give the exact transition automaton on these forms and its extension to every modulus $3^k$, and prove that for integers drawn from a complete residue system modulo $2^j3^k$ the first $j$ residues modulo $3^k$ form \emph{exactly} a Markov chain; its stationary law is Tao's Syracuse random variable, so the form $9n+8$ receives the fraction $22/63$ of all visits. (ii)~We determine the spectrum of this chain completely: its characteristic polynomial is $(x-1)\prod_{i=0}^{k-1}\bigl(x^{2\cdot 3^i}-4^{-3^i}\bigr)$, so every non-trivial eigenvalue has modulus exactly $\tfrac12$ for every $k$. We prove that the mixing time lies between $k\log_2\frac32$ and $2k+O(\sqrt{k})$, show that the stationary mass of the class $-1 \bmod 3^k$ is at least $\tfrac43 2^{-k}$, and present numerical evidence for a cutoff at $k\log_2 3$. (iii)~We explain the strong bias observed in the residue class of the first power of~$2$ met by a trajectory: this class is determined by the last odd term $q_k=(4^k-1)/3$ through $k \bmod 3$; three of the six classes contain only powers of~$2$, a fourth has $O(\log^2N)$ elements up to $N$, and the remaining two are unions of predecessor sets of the $q_k$. Computations for $n\le 4\times10^{9}$ give the frequencies $0.9761$ and $0.0239$. (iv)~We formulate an \emph{octave law}---the expected number of orbit points in $[y,2y)$ lying in a class $r \bmod 3^k$ tends to $\frac{2\log 2}{\log(4/3)}$ times the Syracuse probability of $r$---and support it by exhaustive and Monte Carlo computation.
\end{abstract}



\section{Introduction}\label{sec:intro}

The Collatz map $\Cmap:\N\to\N$ sends $n$ to $n/2$ if $n$ is even and to $3n+1$ if $n$ is odd; the Collatz conjecture asserts that every trajectory reaches~$1$. We refer to Lagarias \cite{lagarias1985,lagarias2010} and Chamberland \cite{chamberland2003} for surveys. The conjecture has been verified for all $n<2^{68}$ \cite{barina2020}. Terras \cite{terras1976} and Everett \cite{everett1977} showed that almost every integer has a finite stopping time, Korec \cite{korec1994} that almost every trajectory falls below $n^{\theta}$ for each $\theta>\log 3/\log 4$, and Tao \cite{tao2022} that almost every trajectory (in logarithmic density) attains values below any function tending to infinity. On the side of predecessor sets, Applegate and Lagarias \cite{applegate1995a,applegate1995b} and Krasikov and Lagarias \cite{krasikov2003} obtained lower bounds for the number of integers below $x$ that reach~$1$, the latter of order $x^{0.84}$, using Krasikov's systems of difference inequalities indexed by residues modulo $3^k$; Wirsching \cite{wirsching1998,wirsching2003} studied predecessor densities through $3$-adic averages.

This paper starts from an elementary observation. Since $2$ is a primitive root modulo~$9$, the powers of $2$ run through the six residues $1,2,4,8,7,5$ modulo~$9$, and once a trajectory has made one odd step it never again meets a multiple of~$3$. Every later term therefore has one of the six \emph{recurrent forms}
\[
 \fa:\ 9n+8,\quad \fb:\ 9n+4,\quad \fc:\ 9n+2,\quad \fd:\ 9n+1,\quad \fe:\ 9n+5,\quad \ff:\ 9n+7 ,
\]
which are the residues of $2^3,2^2,2^1,2^0,2^5,2^4$. Two natural statistics can be attached to them: how often a trajectory \emph{visits} each form, and in which form it \emph{terminates}, i.e.\ which form contains the first power of $2$ it meets. Numerically both are far from uniform, the second spectacularly so: about $97.6\%$ of all integers terminate in form~$\fa$. The purpose of the paper is to explain both phenomena, to place them in the framework of Tao's Syracuse random variables, and to isolate what is provable from what is conjectural.

\subsection*{Results}

\begin{enumerate}[label=(\roman*)]
 \item \emph{Exact automaton and exact Markov property} (Sections~\ref{sec:automaton}--\ref{sec:markov}). \Cref{thm:fsm} lists the transitions between the six forms together with the affine maps acting on the parameter~$n$. For general $k$ the transitions are $r\mapsto r/2$ and $r\mapsto(3r+1)/2$ on $\Z/3^k\Z$. \Cref{thm:exact} shows that, for $n$ ranging over a complete residue system modulo $2^j3^k$, the vector $(n,T(n),\dots,T^j(n))$ reduced modulo $3^k$ is distributed exactly as $j$ steps of the Markov chain $P_k$ that uses the two maps with probability $\frac12$ each, started from the uniform law. The stationary law $\pi_k$ of $P_k$ is the law of Tao's Syracuse random variable (\Cref{prop:stationary}); in particular the six forms $\fa,\fc,\fb,\fd,\fe,\ff$ receive the stationary masses $\frac{22}{63},\frac{16}{63},\frac{11}{63},\frac{8}{63},\frac{4}{63},\frac{2}{63}$.
 \item \emph{Spectrum and mixing} (Section~\ref{sec:spectrum}). \Cref{thm:spectrum} computes the characteristic polynomial of $P_k$; all eigenvalues other than $1$ have modulus $\frac12$, independently of~$k$. The chain is far from normal, however, and its total-variation mixing time grows linearly in $k$: \Cref{thm:mixing} gives the bounds $k\log_2\frac32-O(1)\le t_{\mathrm{mix}}\le 2k+O(\sqrt k)$, the lower bound improving to $k\log_2 3-o(k)$ if the $\ell^2$-norms of the densities $3^k\pi_k$ grow subexponentially. Exact computation for $k\le 12$ supports a cutoff at $k\log_2 3$ (\Cref{conj:cutoff}). We also show $\pi_k(-1)\ge\frac43\,2^{-k}$, so that the densities $3^k\pi_k$ are unbounded; the class $-1 \bmod 3^k$, which generalises form~$\fa$, is the mode of $\pi_k$ for all $k\le 12$.
 \item \emph{Terminal forms} (Section~\ref{sec:terminal}). \Cref{thm:terminal} describes the six termination classes exactly. The first power of $2$ met by a trajectory is $2^{2k-1}$, where $q_k=(4^k-1)/3$ is its last odd term, and the form is decided by $k\bmod 3$. Forms $\fb,\fd,\ff$ are reached only by powers of~$2$; form $\fe$ only by powers of $2$ and by the numbers $2^iq_{3j}$; forms $\fa$ and $\fc$ collect the predecessor sets of $q_k$ with $k\equiv2$ and $k\equiv 1\pmod 3$ respectively. Thus the ``bias'' is the statement that the branch of the Collatz tree hanging from $5$ carries about $93.8\%$ of the integers.
 \item \emph{Predecessor densities and the octave law} (Section~\ref{sec:octave}). We tabulate the empirical densities $\delta_k$ of the predecessor sets of $q_k$ up to $4\times10^9$ and propose \Cref{conj:octave}, which predicts the expected number of orbit points per octave in each class modulo $3^k$. It is confirmed to three digits by a Monte Carlo experiment at heights $2^{27}$--$2^{60}$, and it accounts for the irregular sizes of the $\delta_k$.
\end{enumerate}

\subsection*{Relation to earlier work}
The Markov description of residues modulo $3^k$ is implicit in much of the literature: the transition structure modulo $3^k$ is the forward counterpart of the backward recursions underlying the Krasikov inequalities \cite{applegate1995b,krasikov2003}, and the limiting law was introduced and studied in depth by Tao \cite{tao2022}, whose values for the modulus~$9$ we recover. The elementary facts of Sections~\ref{sec:prelim} and~\ref{sec:automaton} are classical or routine. To our knowledge the explicit characteristic polynomial (\Cref{thm:spectrum}), the mixing-time bounds and cutoff data, the exact description of the termination classes (\Cref{thm:terminal}) and the formulation and numerical test of the octave law have not appeared before; we have not, however, been able to exclude that some of them are known to specialists. The chain $P_k$ is an affine random walk on a cyclic group, a relative of the walks $x\mapsto ax+b$ of Chung, Diaconis and Graham \cite{chung1987}.

\section{Preliminaries}\label{sec:prelim}

Throughout, $\N=\{1,2,3,\dots\}$ and $T:\N\to\N$ is the accelerated map
\[
 T(n)=\begin{cases} n/2,& n \text{ even},\\ (3n+1)/2,& n\text{ odd}.\end{cases}
\]
Following the notation of the six forms we write $E$ for a halving step and $OE$ for an odd step $n\mapsto(3n+1)/2$ (the operation $O:n\mapsto 3n+1$ followed by $E$). The \emph{orbit} of $n$ is the set $\mathcal{O}(n)=\{T^i(n):i\ge0\}$, and the \emph{parity vector} of length $j$ is $(\varepsilon_0(n),\dots,\varepsilon_{j-1}(n))$ with $\varepsilon_i(n)=T^i(n)\bmod 2$.

\begin{lemma}[Terras \cite{terras1976}, Everett \cite{everett1977}]\label{lem:terras}
 Let $j\ge1$. The parity vector of length $j$ depends only on $n \bmod 2^j$, and $n\bmod 2^j\mapsto(\varepsilon_0,\dots,\varepsilon_{j-1})$ is a bijection from $\Z/2^j\Z$ onto $\{0,1\}^j$. Moreover
 \begin{equation}\label{eq:affine}
  T^j(n)=\frac{3^{m}n+c}{2^j},\qquad m=\varepsilon_0+\dots+\varepsilon_{j-1},
 \end{equation}
 where the integer $c\ge0$ depends only on the parity vector.
\end{lemma}

\begin{lemma}[Runs of odd steps]\label{lem:runs}
 For all $j,m\ge1$ we have $T^j(2^jm-1)=3^jm-1$.
\end{lemma}
\begin{proof}
 Since $2^jm-1$ is odd, $T(2^jm-1)=(3\cdot2^jm-2)/2=2^{j-1}(3m)-1$; iterate.
\end{proof}

The case $j=2$ reads $T^2(4m-1)=9m-1$: every $n\equiv3\pmod 4$ reaches the form $\fa$ in two steps, at a larger value. The other residues modulo $4$ descend: $T(n)<n$ for even $n$ and $T^2(4m+1)=3m+1<4m+1$ for $m\ge1$. By strong induction one obtains the classical reduction: \emph{the Collatz conjecture holds if and only if the trajectory of every $n\equiv3\pmod4$ reaches~$1$}.

For an integer $a\ge1$ let $\Pred(a)=\{n\in\N: a\in\mathcal{O}(n)\}$ be its \emph{predecessor set}. If $a$ is not periodic under $T$, then
\begin{equation}\label{eq:predsplit}
 \Pred(a)=\{a\}\ \sqcup\ \Pred(2a)\ \sqcup\ \begin{cases}\Pred\bigl(\tfrac{2a-1}{3}\bigr),& a\equiv2\pmod 3,\\ \emptyset,&\text{otherwise,}\end{cases}
\end{equation}
because the only possible $T$-preimages of $a$ are $2a$ and, when it is an odd integer, $(2a-1)/3$; note that $(2a-1)/3$ is automatically odd when it is an integer. In particular $\Pred(a)=\{2^ia:i\ge0\}$ when $3\mid a$ (a multiple of $3$ is never periodic, since all its iterated preimages are even).

\section{The residue automaton}\label{sec:automaton}

Fix $k\ge1$. Since $2$ is invertible modulo $3^k$, we may define two maps of $\Z/3^k\Z$ into itself,
\begin{equation}\label{eq:phis}
 \varphi_0(r)=2^{-1}r,\qquad \varphi_1(r)=2^{-1}(3r+1).
\end{equation}

\begin{lemma}\label{lem:residue}
 For every $x\in\N$, $\;T(x)\bmod 3^k=\varphi_{x\bmod 2}(x\bmod 3^k)$. Consequently
 \begin{equation}\label{eq:composition}
  T^j(x)\bmod 3^k=\varphi_{\varepsilon_{j-1}}\circ\dots\circ\varphi_{\varepsilon_0}\,(x\bmod 3^k),
 \end{equation}
 where $(\varepsilon_0,\dots,\varepsilon_{j-1})$ is the parity vector of $x$. If the parity vector contains $m$ ones, the right-hand side depends on $x \bmod 3^k$ only through $x\bmod 3^{\max(k-m,0)}$.
\end{lemma}
\begin{proof}
 The first assertion is the definition of $T$ read modulo $3^k$, and \eqref{eq:composition} follows by induction. A composition of the $\varphi_\varepsilon$ containing $m$ maps $\varphi_1$ has the form $r\mapsto 2^{-j}(3^mr+c)$, which proves the last assertion.
\end{proof}

Thus the pair (parity, residue modulo $3^k$) of $x$ determines the residue of $T(x)$: this is a deterministic automaton $\mathcal{A}_k$ with state set $\Z/3^k\Z$ and input alphabet $\{0,1\}$. A residue divisible by $3$ is mapped by $\varphi_0$ to a residue divisible by $3$ and by $\varphi_1$ to a unit, and units are mapped to units; hence the units $(\Z/3^k\Z)^\times$ form a closed set of states and all other states are transient. For $k=2$ the closed part consists of the six recurrent forms, and the parity of $x=9n+r$ is that of $n+r$. Writing the transitions in terms of the parameter $n$ gives the following table.

\begin{theorem}[Automaton on the six forms]\label{thm:fsm}
 Let $x=9n+r$ with $n\ge0$ and $r\in\{8,4,2,1,5,7\}$. Then $T(x)=9n'+r'$, where the operation performed, the new form $r'$ and the new parameter $n'$ are given in \Cref{tab:fsm}.
\end{theorem}

\begin{table}[H]
 \centering
 \caption{Transitions of the six recurrent forms under $T$. An entry ``$E\to\fb$, $n/2$'' means that $x$ is even, $T(x)=x/2$ lies in form $\fb$, and its parameter is $n'=n/2$.}
 \label{tab:fsm}
 \begin{tabular}{llll}
  \toprule
  \textbf{Form of $x$} & \textbf{Residue} & \textbf{$n$ even} & \textbf{$n$ odd} \\
  \midrule
  $\fa:\ 9n+8$ & $8\equiv2^3$ & $E\to\fb$,\quad $n/2$ & $OE\to\fa$,\quad $(3n+1)/2$ \\
  $\fb:\ 9n+4$ & $4\equiv2^2$ & $E\to\fc$,\quad $n/2$ & $OE\to\fc$,\quad $(3n+1)/2$ \\
  $\fc:\ 9n+2$ & $2\equiv2^1$ & $E\to\fd$,\quad $n/2$ & $OE\to\fa$,\quad $(3n-1)/2$ \\
  $\fd:\ 9n+1$ & $1\equiv2^0$ & $OE\to\fc$,\quad $3n/2$ & $E\to\fe$,\quad $(n-1)/2$ \\
  $\fe:\ 9n+5$ & $5\equiv2^5$ & $OE\to\fa$,\quad $3n/2$ & $E\to\ff$,\quad $(n-1)/2$ \\
  $\ff:\ 9n+7$ & $7\equiv2^4$ & $OE\to\fc$,\quad $(3n+2)/2$ & $E\to\fa$,\quad $(n-1)/2$ \\
  \bottomrule
 \end{tabular}
\end{table}

\begin{proof}
 The number $x=9n+r$ is even exactly when $n\equiv r\pmod2$, which gives the operation. Each entry is then a one-line verification; for instance, for $\ff$ with $n$ even, $x$ is odd and $T(x)=(27n+22)/2=9\cdot\frac{3n+2}{2}+2$, and for $\ff$ with $n$ odd, $T(x)=(9n+7)/2=9\cdot\frac{n-1}{2}+8$. The remaining ten cases are identical in nature.
\end{proof}

\begin{figure}[ht]
\centering
\begin{tikzpicture}[>=Stealth, st/.style={circle,draw,minimum size=8.5mm,inner sep=0pt,font=\small}]
 \node[st] (a) at (90:2.3) {$\fa$};
 \node[st] (b) at (30:2.3) {$\fb$};
 \node[st] (c) at (-30:2.3) {$\fc$};
 \node[st] (d) at (-90:2.3) {$\fd$};
 \node[st] (e) at (-150:2.3) {$\fe$};
 \node[st] (f) at (150:2.3) {$\ff$};
 \foreach \u/\v in {a/b,b/c,c/d,d/e,e/f,f/a} \draw[->,thick] (\u) -- (\v);
 \draw[->,dashed] (a) edge[loop above,looseness=7] (a);
 \draw[->,dashed] (b) edge[bend left=40] (c);
 \draw[->,dashed] (c) edge[bend left=12] (a);
 \draw[->,dashed] (d) edge[bend right=40] (c);
 \draw[->,dashed] (e) edge[bend left=12] (a);
 \draw[->,dashed] (f) edge[bend left=8] (c);
\end{tikzpicture}
\caption{The automaton $\mathcal{A}_2$ restricted to the six recurrent forms. Solid arrows: halving steps $E$ (multiplication by $2^{-1}\equiv5$ modulo~$9$), which form a hexagon. Dashed arrows: odd steps $OE$, all of which lead into $\fa$ or $\fc$.}
\label{fig:fsm}
\end{figure}

\begin{remark}\label{rem:ratios}
 In every row of \Cref{tab:fsm} the parameter is multiplied by $\frac12$ (halving step) or by $\frac32$ (odd step), up to an additive constant of modulus at most~$1$. The automaton by itself therefore carries no contraction information beyond the classical one: with fair parities the mean logarithmic drift per step is $\frac12\log\frac34<0$.
\end{remark}

\section{An exact Markov chain and its stationary law}\label{sec:markov}

Let $P_k$ be the transition matrix on $\Z/3^k\Z$ that moves $r$ to $\varphi_0(r)$ or to $\varphi_1(r)$ with probability $\frac12$ each:
\[
 P_k(r,r')=\tfrac12\mathbf 1[r'=\varphi_0(r)]+\tfrac12\mathbf 1[r'=\varphi_1(r)] .
\]

\begin{theorem}[Exact Markov property]\label{thm:exact}
 Let $j\ge0$, $k\ge1$ and $r_0,\dots,r_j\in\Z/3^k\Z$. The set
 \[
  A(r_0,\dots,r_j)=\{n\in\N:\ T^i(n)\equiv r_i \!\!\pmod{3^k}\ \text{for } 0\le i\le j\}
 \]
 is a union of residue classes modulo $2^j3^k$, and its natural density equals
 \[
  3^{-k}\,\prod_{i=0}^{j-1}P_k(r_i,r_{i+1}).
 \]
 Equivalently: if $n$ is uniformly distributed over a complete residue system modulo $2^j3^k$, then $(T^i(n)\bmod 3^k)_{0\le i\le j}$ has exactly the law of the first $j$ steps of the Markov chain with kernel $P_k$ started from the uniform distribution.
\end{theorem}
\begin{proof}
 By the Chinese remainder theorem, when $n$ is uniform modulo $2^j3^k$ the residues $n\bmod 2^j$ and $n\bmod3^k$ are independent and uniform. By \Cref{lem:terras} the parity vector $(\varepsilon_0,\dots,\varepsilon_{j-1})$ is a bijective function of $n\bmod 2^j$, hence is uniform on $\{0,1\}^j$ and independent of $n\bmod3^k$. By \eqref{eq:composition}, $T^{i}(n)\bmod3^k=\varphi_{\varepsilon_{i-1}}(T^{i-1}(n)\bmod 3^k)$ for $1\le i\le j$. This is precisely the chain $P_k$ driven by independent fair bits and started from a uniform state. Since membership in $A(r_0,\dots,r_j)$ depends only on $n\bmod 2^j3^k$, the density statement follows.
\end{proof}

\begin{proposition}[Stationary law]\label{prop:stationary}
 Let $k\ge1$.
 \begin{enumerate}[label=(\alph*)]
  \item $P_k$ has a unique stationary distribution $\pi_k$, and for every initial distribution $\mu$ and every $t\ge0$,
  \begin{equation}\label{eq:coupling}
   \|\mu P_k^t-\pi_k\|_{\TV}\le\Prob\bigl(\Bin(t,\tfrac12)<k\bigr).
  \end{equation}
  \item $\pi_k$ is the law of
  \begin{equation}\label{eq:syrac}
   X_k=\sum_{i=1}^{k}3^{i-1}\,2^{-(g_1+\dots+g_i)} \pmod{3^k},
  \end{equation}
  where $g_1,g_2,\dots$ are independent with $\Prob(g_i=g)=2^{-g}$, $g\ge1$. That is, $\pi_k$ is the law of the Syracuse random variable $\mathbf{Syrac}(\Z/3^k\Z)$ of Tao \cite{tao2022}.
  \item With $\pi_0\equiv1$ on the trivial group, for every $x\in\Z/3^k\Z$,
  \begin{equation}\label{eq:recursion}
   \pi_k(x)=\sum_{\substack{g\ge1\\ 2^gx\equiv1\ (\mathrm{mod}\ 3)}}2^{-g}\,\pi_{k-1}\!\Bigl(\frac{2^gx-1}{3}\bmod 3^{k-1}\Bigr).
  \end{equation}
  In particular $\pi_k(x)>0$ if and only if $3\nmid x$.
  \item $\displaystyle \tfrac43\,2^{-k}\le\pi_k(-1)\le\max_x\pi_k(x)\le\bigl(\tfrac23\bigr)^k$, and $2^k\pi_k(-1)$ is non-decreasing in $k$.
 \end{enumerate}
\end{proposition}

\begin{proof}
 (a) Run two copies of the chain from arbitrary initial states with the same driving bits. By the last assertion of \Cref{lem:residue}, as soon as $k$ of the bits equal $1$ the two copies coincide, and they coincide forever after. Hence the chain has at most one stationary distribution; one exists because the state space is finite. Taking the second copy stationary gives \eqref{eq:coupling} by the coupling inequality.

 (b) By (a), $\pi_k$ is the limit in law of $\varphi_{\varepsilon_t}\circ\dots\circ\varphi_{\varepsilon_1}(0)$ for independent fair bits $\varepsilon_i$. Reading the bits from the last one backwards, let $g_1-1$ be the number of zeros before the first one, $g_2-1$ the number of zeros between the first and the second one, and so on; these $g_i$ are independent and geometric as stated. On the event that at least $k$ ones occur, the composition is $\varphi_0^{g_1-1}\varphi_1\varphi_0^{g_2-1}\varphi_1\cdots\varphi_0^{g_k-1}\varphi_1$ applied to some residue $y$, and since $\varphi_0^{g-1}\varphi_1(z)=2^{-g}(3z+1)$ this equals $\sum_{i=1}^k3^{i-1}2^{-(g_1+\dots+g_i)}+3^k2^{-(g_1+\dots+g_k)}y\equiv X_k$. Letting $t\to\infty$ proves (b); the expression \eqref{eq:syrac} is Tao's definition up to reversing the order of the $g_i$.

 (c) By \eqref{eq:syrac}, $X_k\equiv2^{-g_1}(1+3X')$, where $X'=\sum_{i=2}^k3^{i-2}2^{-(g_2+\dots+g_i)}$ is independent of $g_1$ and has law $\pi_{k-1}$ modulo $3^{k-1}$. For fixed $g_1=g$, the congruence $2^{-g}(1+3X')\equiv x\pmod{3^k}$ holds if and only if $2^gx\equiv1\pmod3$ and $X'\equiv(2^gx-1)/3\pmod{3^{k-1}}$. This gives \eqref{eq:recursion}. If $3\mid x$ the sum is empty. If $3\nmid x$, then as $g$ runs through its admissible parity class, $2^gx$ runs through all residues modulo $3^k$ that are $\equiv1\pmod 3$ (because $4$ generates the group of such residues), so $(2^gx-1)/3$ runs through all of $\Z/3^{k-1}\Z$ and positivity follows by induction.

 (d) In \eqref{eq:recursion} the admissible $g$ are either all odd or all even, so $\sum 2^{-g}\le\frac23$ and $\max\pi_k\le\frac23\max\pi_{k-1}$. For $x=-1$ the admissible $g$ are odd and $g=1$ contributes $\frac12\pi_{k-1}(-1)$, whence $2^k\pi_k(-1)\ge2^{k-1}\pi_{k-1}(-1)$; and $\pi_1(-1)=\sum_{g\text{ odd}}2^{-g}=\frac23$.
\end{proof}

For $k=1$ the recursion gives $\pi_1=(0,\frac13,\frac23)$ on the residues $(0,1,2)$, and for $k=2$
\begin{equation}\label{eq:pi2}
 \pi_2(\fa)=\tfrac{22}{63},\quad \pi_2(\fc)=\tfrac{16}{63},\quad \pi_2(\fb)=\tfrac{11}{63},\quad \pi_2(\fd)=\tfrac{8}{63},\quad \pi_2(\fe)=\tfrac{4}{63},\quad \pi_2(\ff)=\tfrac{2}{63},
\end{equation}
in agreement with the values listed in \cite[Section~1.3]{tao2022}. These can also be read off \Cref{fig:fsm}: each state sends half of its mass along its solid arrow and half along its dashed arrow, so $\pi(\fb)=\frac12\pi(\fa)$, $\pi(\fd)=\frac12\pi(\fc)$, $\pi(\fe)=\frac12\pi(\fd)$, $\pi(\ff)=\frac12\pi(\fe)$ and $\pi(\fc)=\pi(\fb)+\frac12\pi(\fd)+\frac12\pi(\ff)$.

\begin{corollary}[Visit frequencies]\label{cor:visits}
 For every $j\ge0$, $k\ge1$ and $r\in\Z/3^k\Z$ the set $\{n:T^j(n)\equiv r\pmod{3^k}\}$ has a natural density $\mu_j(r)$, and
 \[
  \sum_{r}\bigl|\mu_j(r)-\pi_k(r)\bigr|\le2\,\Prob\bigl(\Bin(j,\tfrac12)<k\bigr)\le 2\exp\!\Bigl(-\frac{(j-2k)^2}{2j}\Bigr)\quad(j\ge2k).
 \]
 In particular the density of integers whose $j$th iterate has the form $9n+8$ tends to $\frac{22}{63}\approx0.349$ as $j\to\infty$, compared with $\frac{2}{63}\approx0.032$ for the form $9n+7$.
\end{corollary}
\begin{proof}
 Combine \Cref{thm:exact}, \eqref{eq:coupling} and Hoeffding's inequality.
\end{proof}

\begin{remark}[Why form $\fa$ dominates the visits]\label{rem:minusone}
 The class $-1\bmod 3^k$ is the fixed point of the odd-step map $\varphi_1$. By \Cref{lem:runs}, a run of $k$ consecutive odd steps always ends in this class, and it is left only by a halving step. Part~(d) of \Cref{prop:stationary} quantifies the resulting accumulation: the class $-1$ carries at least $\frac43(\frac32)^k$ times the uniform mass $3^{-k}$. For $k=2$ this class is form~$\fa$.
\end{remark}

\section{Spectrum and mixing}\label{sec:spectrum}

\begin{theorem}[Spectrum]\label{thm:spectrum}
 For every $k\ge1$ the characteristic polynomial of $P_k$ is
 \begin{equation}\label{eq:charpoly}
  \det(xI-P_k)=(x-1)\prod_{i=0}^{k-1}\Bigl(x^{2\cdot3^{i}}-4^{-3^{i}}\Bigr).
 \end{equation}
 The characteristic polynomial of the restriction of $P_k$ to the closed class $(\Z/3^k\Z)^\times$ is
 \[
  (x-1)\,\frac{x^{2\cdot3^{k-1}}-4^{-3^{k-1}}}{x-\frac12}.
 \]
 In particular every eigenvalue of $P_k$ other than $1$ has modulus exactly $\frac12$.
\end{theorem}

\begin{proof}
 Let $e(u)=\exp(2\pi iu)$ and, for $\xi\in\Z/3^k\Z$, let $\chi_\xi(r)=e(\xi r/3^k)$. Acting on functions, $(P_kf)(r)=\frac12f(\varphi_0(r))+\frac12f(\varphi_1(r))$, so by \eqref{eq:phis}
 \begin{equation}\label{eq:fourier}
  P_k\chi_\xi=\tfrac12\,\chi_{2^{-1}\xi}+\tfrac12\,e\bigl(2^{-1}\xi/3^k\bigr)\,\chi_{3\cdot2^{-1}\xi}.
 \end{equation}
 Let $v(\xi)\in\{0,\dots,k\}$ be the $3$-adic valuation of $\xi$ (with $v(0)=k$) and $\Xi_v=\{\xi:v(\xi)=v\}$. In \eqref{eq:fourier} the first frequency has valuation $v(\xi)$ and the second has valuation $\min(v(\xi)+1,k)$. Ordering the basis $(\chi_\xi)$ by increasing valuation, the matrix of $P_k$ is therefore block triangular, and its diagonal block on $\Xi_v$ is $\frac12\Pi_v$, where $\Pi_v$ is the permutation $\xi\mapsto2^{-1}\xi$ of $\Xi_v$; for $v=k$ the block is the scalar $\frac12+\frac12=1$. For $v<k$, $\Xi_v=3^v(\Z/3^{k-v}\Z)^\times$ has $2\cdot3^{k-v-1}$ elements, and since $2$ is a primitive root modulo every power of $3$, $\Pi_v$ is a single cycle of that length $L_v$. The characteristic polynomial of $\frac12\Pi_v$ is $x^{L_v}-2^{-L_v}$, and \eqref{eq:charpoly} follows on putting $i=k-1-v$.

 For the restriction, order the states as units first, then multiples of $3$. Units are mapped to units, a multiple of $3$ is mapped by $\varphi_0$ to a multiple of $3$ and by $\varphi_1$ to a unit. Hence $P_k$ is block triangular with diagonal blocks $P_k|_{\text{units}}$ and $\frac12\Pi'$, where $\Pi'$ is the permutation $r\mapsto2^{-1}r$ of $3\Z/3^k\Z$. Its cycles are $\{0\}$ and, for $1\le w\le k-1$, one cycle of length $2\cdot3^{k-w-1}$, so $\det(xI-\frac12\Pi')=(x-\frac12)\prod_{i=0}^{k-2}(x^{2\cdot3^i}-4^{-3^i})$. Dividing \eqref{eq:charpoly} by this polynomial gives the claim.
\end{proof}

\begin{remark}
 For $k=2$ the six-state chain of \Cref{fig:fsm} has characteristic polynomial $(x-1)(x^6-\frac1{64})/(x-\frac12)$. A spectral gap of $\frac12$ that does not depend on $k$ might suggest mixing in $O(1)$ steps, but $P_k$ is highly non-normal (\eqref{eq:fourier} is a two-term recursion that only ever raises the valuation), and the true mixing time is of order $k$, as we show next.
\end{remark}

Let
\[
 d_k(t)=\max_{r\in\Z/3^k\Z}\bigl\|P_k^t(r,\cdot)-\pi_k\bigr\|_{\TV},\qquad t_{\mathrm{mix}}^{(k)}(\varepsilon)=\min\{t:d_k(t)\le\varepsilon\},
\]
and let $c_k=3^k\sum_x\pi_k(x)^2$ be the squared $\ell^2$-norm of the density $3^k\pi_k$ with respect to the uniform probability measure.

\begin{theorem}[Mixing bounds]\label{thm:mixing}
 For all $k\ge1$ and $t\ge0$,
 \begin{equation}\label{eq:mixbounds}
  1-\bigl(2^t\,c_k\,3^{-k}\bigr)^{1/2}\ \le\ d_k(t)\ \le\ \Prob\bigl(\Bin(t,\tfrac12)<k\bigr),
 \end{equation}
 and $c_k\le2^k$. Consequently, for each fixed $\varepsilon\in(0,1)$,
 \[
  k\log_2\tfrac32-2\log_2\tfrac1{1-\varepsilon}\ \le\ t^{(k)}_{\mathrm{mix}}(\varepsilon)\ \le\ 2k+O_\varepsilon(\sqrt k),
 \]
 and if $c_k=2^{o(k)}$ then $t^{(k)}_{\mathrm{mix}}(\varepsilon)\ge k\log_23-o(k)$.
\end{theorem}
\begin{proof}
 The upper bound is \eqref{eq:coupling}; the binomial tail is at most $\varepsilon$ once $t\ge2k+\sqrt{2t\log(1/\varepsilon)}$ by Hoeffding's inequality. For the lower bound, the support $S$ of $P_k^t(r,\cdot)$ has at most $2^t$ elements, so by the Cauchy--Schwarz inequality
 \[
  d_k(t)\ge P_k^t(r,S)-\pi_k(S)=1-\pi_k(S)\ge1-\Bigl(|S|\sum_x\pi_k(x)^2\Bigr)^{1/2}\ge 1-\bigl(2^tc_k3^{-k}\bigr)^{1/2}.
 \]
 Finally $c_k\le3^k\max_x\pi_k(x)\le2^k$ by \Cref{prop:stationary}(d).
\end{proof}

\Cref{tab:mixing} reports exact computations of the stationary laws and of the worst-case distance. The times at which $d_k$ crosses $\frac12$, $\frac14$ and $0.05$ grow at a rate consistent with $\log_23\approx1.585$ per unit of $k$ and stay well below the coupling bound, and the profiles $t\mapsto d_k(t)$ sharpen as $k$ grows when plotted against $t/k$ (\Cref{fig:cutoff}). The constants $c_k$ grow only linearly, $c_k\approx0.47\,k+1.2$, so the conditional lower bound of \Cref{thm:mixing} applies numerically.

\begin{table}[H]
 \centering
 \caption{Mixing data for $P_k$. Columns 2--4: first $t$ with $d_k(t)\le\varepsilon$. For $k\le8$ the maximum defining $d_k$ is taken over all $3^k$ initial states; for $9\le k\le12$ it is taken over nine initial states ($1,2,4,5,7,-1,-2,\pm\frac12$), so that those entries are lower bounds. Column 5: smallest $t$ with $\Prob(\Bin(t,\frac12)<k)\le\frac14$.}
 \label{tab:mixing}
 \begin{tabular}{rrrrrrrr}
  \toprule
  $k$ & $\varepsilon=\frac12$ & $\varepsilon=\frac14$ & $\varepsilon=0.05$ & coupling & $k\log_23$ & $2^k\pi_k(-1)$ & $c_k$ \\
  \midrule
  1 & 1 & 2 & 5 & 2 & 1.58 & 1.3333 & 1.6667 \\
  2 & 3 & 4 & 5 & 5 & 3.17 & 1.3968 & 2.1429 \\
  3 & 4 & 5 & 8 & 7 & 4.75 & 1.4212 & 2.6044 \\
  4 & 6 & 7 & 9 & 10 & 6.34 & 1.4358 & 3.0686 \\
  5 & 7 & 8 & 11 & 12 & 7.92 & 1.4456 & 3.5342 \\
  6 & 9 & 10 & 13 & 14 & 9.51 & 1.4509 & 4.0003 \\
  7 & 10 & 12 & 15 & 16 & 11.09 & 1.4543 & 4.4658 \\
  8 & 12 & 14 & 17 & 18 & 12.68 & 1.4567 & 4.9317 \\
  9 & 13 & 15 & 18 & 21 & 14.26 & 1.4587 & 5.3984 \\
  10 & 15 & 17 & 20 & 23 & 15.85 & 1.4598 & 5.8658 \\
  11 & 17 & 19 & 22 & 25 & 17.43 & 1.4606 & 6.3340 \\
  12 & 18 & 20 & 24 & 27 & 19.02 & 1.4613 & 6.8029 \\
  \bottomrule
 \end{tabular}
\end{table}

\begin{figure}[H]
\centering
\begin{tikzpicture}
\begin{axis}[width=10.5cm,height=5.6cm,xlabel={$t/k$},ylabel={$d_k(t)$},xmin=0,xmax=3.2,ymin=0,ymax=1.02,
 legend style={draw=none,font=\small},legend pos=north east,grid=both,grid style={gray!20},tick label style={font=\small}]
\addplot[mark=o,mark size=1.4pt,gray] coordinates {(0.2500,0.9964) (0.5000,0.9831) (0.7500,0.9275) (1.0000,0.7693) (1.2500,0.5197) (1.5000,0.2937) (1.7500,0.1632) (2.0000,0.0905) (2.2500,0.0486) (2.5000,0.0261) (2.7500,0.0142) (3.0000,0.0071)};
\addplot[mark=square*,mark size=1.2pt,black!60] coordinates {(0.1250,1.0000) (0.2500,0.9999) (0.3750,0.9994) (0.5000,0.9979) (0.6250,0.9937) (0.7500,0.9825) (0.8750,0.9614) (1.0000,0.9175) (1.1250,0.8329) (1.2500,0.7089) (1.3750,0.5524) (1.5000,0.3910) (1.6250,0.2539) (1.7500,0.1566) (1.8750,0.0940) (2.0000,0.0563) (2.1250,0.0340) (2.2500,0.0201) (2.3750,0.0118) (2.5000,0.0069) (2.6250,0.0039) (2.7500,0.0022) (2.8750,0.0013) (3.0000,0.0007) (3.1250,0.0004)};
\addplot[mark=*,mark size=1.2pt,black,thick] coordinates {(0.0833,1.0000) (0.1667,1.0000) (0.2500,1.0000) (0.3333,1.0000) (0.4167,0.9996) (0.5000,0.9993) (0.5833,0.9986) (0.6667,0.9973) (0.7500,0.9945) (0.8333,0.9889) (0.9167,0.9813) (1.0000,0.9694) (1.0833,0.9475) (1.1667,0.9076) (1.2500,0.8413) (1.3333,0.7416) (1.4167,0.6102) (1.5000,0.4631) (1.5833,0.3252) (1.6667,0.2160) (1.7500,0.1390) (1.8333,0.0875) (1.9167,0.0545) (2.0000,0.0339) (2.0833,0.0210) (2.1667,0.0129) (2.2500,0.0078) (2.3333,0.0047) (2.4167,0.0028) (2.5000,0.0017) (2.5833,0.0010) (2.6667,0.0006) (2.7500,0.0004) (2.8333,0.0002) (2.9167,0.0001) (3.0000,0.0001) (3.0833,0.0000) (3.1667,0.0000)};
\draw[dashed] (axis cs:1.585,0) -- (axis cs:1.585,1.02);
\legend{$k=4$,$k=8$,$k=12$}
\end{axis}
\end{tikzpicture}
\caption{Worst-case total-variation distance to stationarity against $t/k$ (for $k=12$, maximum over the nine initial states of \Cref{tab:mixing}). The dashed line marks $t/k=\log_23$.}
\label{fig:cutoff}
\end{figure}

\begin{conjecture}[Cutoff]\label{conj:cutoff}
 The chains $P_k$ exhibit a total-variation cutoff at $k\log_23$: for every $\eta>0$, $d_k\bigl((1-\eta)k\log_23\bigr)\to1$ and $d_k\bigl((1+\eta)k\log_23\bigr)\to0$ as $k\to\infty$.
\end{conjecture}

The heuristic is entropic: $t$ fair bits can select at most $2^t$ residues, and $\pi_k$ is spread over about $3^k$ of them. The lower half of the conjecture would follow from $c_k=2^{o(k)}$; the upper half is a statement of equidistribution at the finest scale, of the same nature as the fine-scale mixing estimate that is the technical heart of \cite{tao2022}, and appears to be difficult.

\begin{remark}[Unbounded densities]\label{rem:unbounded}
 By \Cref{prop:stationary}(d), $\max_x3^k\pi_k(x)\ge\frac43(\frac32)^k$: the densities of the Syracuse laws with respect to the uniform measure are unbounded, with a spike at $-1$ (and, by stationarity, smaller spikes at its images under compositions of $\varphi_0,\varphi_1$). The seventh column of \Cref{tab:mixing} suggests that $2^k\pi_k(-1)$ converges, to a limit near $1.46$; for every $k\le12$ the maximum of $\pi_k$ is attained at $-1$. The slow, apparently linear growth of $c_k$ raises the question whether the limiting density is square integrable.
\end{remark}

\section{Terminal forms}\label{sec:terminal}

We now turn to the second statistic. If the orbit of $n$ contains a power of $2$, let $\tau(n)$ be the first power of $2$ in the sequence $n,T(n),T^2(n),\dots$, and for a form $r\in\{\fa,\dots,\ff\}$ let
\[
 N_r=\{n\in\N:\ \tau(n)\text{ exists and has form } r\},\qquad \Freq_r(N)=\frac{|N_r\cap[1,N]|}{N}.
\]
Let $q_k=(4^k-1)/3=1+4+\dots+4^{k-1}$, so that $q_1=1$, $q_2=5$, $q_3=21$, $q_4=85$, $q_5=341$.

\begin{lemma}\label{lem:lastodd}
 Suppose $\tau(n)$ exists and $n$ is not a power of $2$. Then the term preceding $\tau(n)$ in the trajectory of $n$ is $q_k$ for some $k\ge2$, and $\tau(n)=2^{2k-1}$. Moreover $q_k\equiv k\pmod3$.
\end{lemma}
\begin{proof}
 Let $m$ be the term preceding $\tau(n)=2^s$. If $m$ were even we would have $m=2^{s+1}$, contradicting the minimality in the definition of $\tau(n)$. Hence $m$ is odd, $(3m+1)/2=2^s$, and $m=(2^{s+1}-1)/3$ is an integer only for $s+1=2k$ even; $m$ is not a power of $2$, so $k\ge2$. Finally $q_k=\sum_{i<k}4^i\equiv k\pmod 3$.
\end{proof}

\begin{theorem}[Structure of the termination classes]\label{thm:terminal}
 The termination classes are as follows.
 \begin{enumerate}[label=(\roman*)]
  \item $N_{\fb}=\{2^{6j+2}:j\ge0\}$, $N_{\fd}=\{2^{6j}:j\ge0\}$ and $N_{\ff}=\{2^{6j+4}:j\ge0\}$.
  \item $N_{\fe}=\{2^{6j+5}:j\ge0\}\cup\{2^iq_{3j}:i\ge0,\ j\ge1\}$, and consequently
  \[
   |N_{\fe}\cap[1,N]|=\#\{j\ge0:2^{6j+5}\le N\}+\sum_{j\ge1,\ q_{3j}\le N}\Bigl(\floor{\log_2(N/q_{3j})}+1\Bigr)=O(\log^2N).
  \]
  \item $N_{\fa}=\{2^{6j+3}:j\ge0\}\cup\bigcup_{k\equiv2\ (3)}\Pred(q_k)$ and $N_{\fc}=\{2^{6j+1}:j\ge0\}\cup\bigcup_{k\equiv1\ (3),\,k\ge4}\Pred(q_k)$, the unions being disjoint.
  \item Let $\rho(N)$ be the proportion of $n\le N$ whose trajectory reaches $1$. Then
  \[
   \Freq_{\fa}(N)+\Freq_{\fc}(N)=\rho(N)-O\!\Bigl(\frac{\log^2N}{N}\Bigr).
  \]
 \end{enumerate}
\end{theorem}
\begin{proof}
 If $n=2^s$ then $\tau(n)=n$, whose form is determined by $s\bmod6$. Otherwise \Cref{lem:lastodd} applies: $\tau(n)=2^{2k-1}$ with $n\in\Pred(q_k)$, $k\ge2$, and conversely, for $k\ge2$, no element of $\Pred(q_k)$ is a power of $2$ (the orbit of a power of $2$ consists of powers of $2$), so every $n\in\Pred(q_k)$ has $\tau(n)=T(q_k)=2^{2k-1}$. Now $2k-1\bmod 6$ equals $3,5,1$ for $k\equiv2,0,1\pmod3$, i.e.\ $2^{2k-1}$ has form $\fa,\fe,\fc$ respectively. This proves (i) and (iii); the sets $\Pred(q_k)$ are pairwise disjoint because the term preceding $\tau(n)$ is unique. For (ii), $k\equiv0\pmod3$ gives $3\mid q_k$ by \Cref{lem:lastodd}, hence $\Pred(q_k)=\{2^iq_k\}$ by \eqref{eq:predsplit}; counting gives the formula. Part (iv) follows since $n$ reaches $1$ if and only if $\tau(n)$ exists.
\end{proof}

\begin{remark}[Dependence on the choice of map]
 If the original map $\Cmap$ is used instead of $T$, the first power of $2$ met by a non-power of $2$ is $3q_k+1=2^{2k}$ rather than $2^{2k-1}$, and the roles of the forms are rotated: $k\equiv2,1,0\pmod 3$ now correspond to $\ff,\fb,\fd$. The partition of the integers according to $k\bmod3$ is the same; only its labelling by forms changes.
\end{remark}

\begin{remark}[What is and is not known about the limits]\label{rem:limits}
 By (iv), $\Freq_{\fa}(N)+\Freq_{\fc}(N)\to1$ holds if and only if the set of integers whose trajectory reaches $1$ has natural density~$1$. This is open: the results of \cite{korec1994,tao2022} concern the minimum of an orbit, not its eventual fate, and the best unconditional lower bound for the number of $n\le N$ reaching $1$ is $N^{0.84}$ \cite{krasikov2003}. For $N<2^{68}$ one has $\rho(N)=1$ by \cite{barina2020}. The existence of the individual limits of $\Freq_{\fa}$ and $\Freq_{\fc}$ is likewise open, since it concerns the densities of the sets $\Pred(q_k)$.
\end{remark}

\Cref{tab:freq} gives the exact counts for $N\le4\times10^9$; the counts for $\fb,\fd,\fe,\ff$ agree with the formulas of \Cref{thm:terminal}. \Cref{tab:delta} splits the classes $\fa$ and $\fc$ according to the last odd term. The branch at $q_2=5$ alone accounts for $93.77\%$ of all integers, the branch at $q_5=341$ for $3.80\%$ and the branch at $q_4=85$ for $2.38\%$.

\begin{table}[H]
 \centering
 \caption{Termination classes for $n\le N$ (first power of $2$ in the $T$-trajectory).}
 \label{tab:freq}
 \begin{tabular}{rcccccc}
  \toprule
  $N$ & $\Freq_{\fa}(N)$ & $\Freq_{\fc}(N)$ & $|N_{\fb}|$ & $|N_{\fd}|$ & $|N_{\fe}|$ & $|N_{\ff}|$ \\
  \midrule
  $10$ & 0.7000 & 0.1000 & 1 & 1 & 0 & 0 \\
  $10^{2}$ & 0.8900 & 0.0300 & 1 & 2 & 4 & 1 \\
  $10^{3}$ & 0.9590 & 0.0290 & 2 & 2 & 7 & 1 \\
  $10^{4}$ & 0.97390 & 0.02400 & 2 & 3 & 14 & 2 \\
  $10^{5}$ & 0.97478 & 0.02490 & 3 & 3 & 23 & 3 \\
  $10^{6}$ & 0.976082 & 0.023875 & 3 & 4 & 33 & 3 \\
  $10^{7}$ & 0.976114 & 0.023881 & 4 & 4 & 44 & 4 \\
  $10^{8}$ & 0.976161 & 0.023838 & 5 & 5 & 60 & 4 \\
  $10^{9}$ & 0.976131 & 0.023869 & 5 & 5 & 75 & 5 \\
  $2\times10^{9}$ & 0.976129 & 0.023871 & 5 & 6 & 80 & 5 \\
  $4\times10^{9}$ & 0.976125 & 0.023875 & 5 & 6 & 85 & 5 \\
  \bottomrule
 \end{tabular}
\end{table}

\begin{table}[H]
 \centering
 \caption{Empirical densities $\delta_k(N)=|\Pred(q_k)\cap[1,N]|/N$ of the terminal branches, the normalised values $q_k\delta_k/\kappa$ with $\kappa=2/\log\frac43$ at $N=4\times10^9$, and the Syracuse densities $9\pi_2$ and $3^6\pi_6$ at $q_k$ (see Section~\ref{sec:octave}).}
 \label{tab:delta}
 \begin{tabular}{rrcccccc}
  \toprule
  $k$ & $q_k$ & form & $q_k\bmod9$ & $\delta_k(10^9)$ & $\delta_k(4\times10^9)$ & $q_k\delta_k/\kappa$ & $9\pi_2,\ 3^6\pi_6$ \\
  \midrule
  2 & 5 & $\mathfrak{a}$ & 5 & $9.377\times10^{-1}$ & $9.377\times10^{-1}$ & 0.674 & 0.571,\ 0.702 \\
  4 & 85 & $\mathfrak{c}$ & 4 & $2.375\times10^{-2}$ & $2.376\times10^{-2}$ & 0.291 & 1.571,\ 0.560 \\
  5 & 341 & $\mathfrak{a}$ & 8 & $3.795\times10^{-2}$ & $3.796\times10^{-2}$ & 1.862 & 3.143,\ 1.868 \\
  7 & 5461 & $\mathfrak{c}$ & 7 & $8.077\times10^{-5}$ & $8.009\times10^{-5}$ & 0.063 & 0.286,\ 0.150 \\
  8 & 21845 & $\mathfrak{a}$ & 2 & $4.860\times10^{-4}$ & $4.853\times10^{-4}$ & 1.525 & 2.286,\ 2.176 \\
  10 & 349525 & $\mathfrak{c}$ & 1 & $3.253\times10^{-5}$ & $3.232\times10^{-5}$ & 1.625 & 1.143,\ 3.569 \\
  11 & 1398101 & $\mathfrak{a}$ & 5 & $2.106\times10^{-6}$ & $2.118\times10^{-6}$ & 0.426 & 0.571,\ 0.426 \\
  \bottomrule
 \end{tabular}
\end{table}

\section{Predecessor densities and the octave law}\label{sec:octave}

\Cref{thm:terminal} reduces the termination statistics to the sizes of the predecessor sets $\Pred(q_k)$. \Cref{tab:delta} shows that these sizes are far from monotone in $k$: for instance $\Pred(q_8)$ is six times larger than $\Pred(q_7)$ although $q_8\approx4q_7$. We now describe a heuristic law that accounts for this and connects the two halves of the paper.

For $y\ge1$, $k\ge1$, $r\in\Z/3^k\Z$ and $N\ge1$ let
\[
 V_{k,r}(y;N)=\frac1N\sum_{n\le N}\#\bigl\{a\in\mathcal{O}(n):\ y\le a<2y,\ a\equiv r\!\!\pmod{3^k}\bigr\}
 =\sum_{\substack{y\le a<2y\\ a\equiv r\ (3^k)}}\frac{|\Pred(a)\cap[1,N]|}{N}
\]
be the average number of orbit points in the octave $[y,2y)$ that lie in the class $r$.

\begin{conjecture}[Octave law]\label{conj:octave}
 For every $k\ge1$ and $r\in\Z/3^k\Z$,
 \[
  \lim_{y\to\infty}\ \liminf_{N\to\infty}V_{k,r}(y;N)=\lim_{y\to\infty}\ \limsup_{N\to\infty}V_{k,r}(y;N)=\frac{2\log2}{\log(4/3)}\ \pi_k(r).
 \]
\end{conjecture}

\subsection*{Heuristic derivation}
By \Cref{lem:terras} and \Cref{thm:exact}, the initial segment of the trajectory of a random integer is driven by independent fair bits: each bit multiplies the current value by $\frac12$ or (essentially) $\frac32$ and simultaneously moves the residue modulo $3^k$ by $\varphi_0$ or $\varphi_1$. The pair (residue, logarithm of the value) is thus a Markov additive process whose additive part has increments $-\log2$ and $\log\frac32$, with mean $-\frac12\log\frac43$ under the stationary law and non-lattice support. If this description remained valid down to height~$y$, the Markov renewal theorem \cite[Ch.~VII]{asmussen2003} would give, for the expected number of visits to a window of logarithmic length $\ell$ in residue class $r$, the value $\ell\,\pi_k(r)/(\frac12\log\frac43)$; with $\ell=\log2$ this is the conjecture. Summing over $r$, the predicted number of orbit points per octave is $2\log 2/\log\frac43\approx4.8188$, which is the classical prediction $\kappa\log n$, $\kappa=2/\log\frac43\approx6.952$, for the total stopping time \cite{lagarias1992,borovkov2000,kontorovich2010}. The rigorous description is only available for the first $\log_2N$ steps of the trajectories of $n\le N$, during which a typical trajectory descends only to height about $N^{\log3/\log4}$ \cite{korec1994}; the conjecture asserts that the same statistics persist at every sufficiently large fixed height.

\subsection*{Numerical evidence}
We tested \Cref{conj:octave} in two ways. First, exhaustively: for every $n\le N=4\times10^9$ we recorded the first orbit point below $2^{15}$, from which $|\Pred(a)\cap[1,N]|$ is obtained for every $a<2^{15}$. Second, by simulation at larger heights: we followed the complete trajectories of $6\times10^4$ random $90$-bit integers and counted orbit points in the octaves $[2^j,2^{j+1})$, $16\le j<60$. The results are given in \Cref{tab:octave}. At heights $2^{27}$ and above the number of points per octave agrees with $4.8188$ to within $0.004$, and the distribution over the six forms agrees with \eqref{eq:pi2} to within $6\times10^{-4}$ in total variation (and to within $9\times 10^{-4}$ modulo $27$), which is the order of the sampling error. At small heights the statistics are visibly lumpy and converge as $y$ grows.

\begin{table}[H]
 \centering
 \caption{Orbit points per octave and their distribution over the six forms. Top: exhaustive, $n\le4\times10^9$, octave $[2^j,2^{j+1})$. Bottom: $6\times10^4$ random $90$-bit starting values, averaged over the indicated range of octaves. Last column: total-variation distance between the observed distribution and $\pi_2$.}
 \label{tab:octave}
 \begin{tabular}{lcccccccc}
  \toprule
  octaves & points & $\fa$ & $\fc$ & $\fb$ & $\fd$ & $\fe$ & $\ff$ & TV \\
  \midrule
  $j=10$ & 5.610 & 0.3843 & 0.2465 & 0.1948 & 0.0877 & 0.0562 & 0.0305 & 0.0553 \\
  $j=11$ & 6.017 & 0.4739 & 0.1720 & 0.1761 & 0.1014 & 0.0529 & 0.0236 & 0.1262 \\
  $j=12$ & 4.736 & 0.3811 & 0.2343 & 0.1544 & 0.1338 & 0.0633 & 0.0331 & 0.0401 \\
  $j=13$ & 4.655 & 0.3339 & 0.2631 & 0.1761 & 0.1257 & 0.0686 & 0.0326 & 0.0166 \\
  $j=14$ & 4.722 & 0.3432 & 0.2518 & 0.1808 & 0.1279 & 0.0643 & 0.0319 & 0.0082 \\
  \midrule
  $16\le j<27$ & 4.773 & 0.3475 & 0.2545 & 0.1744 & 0.1279 & 0.0637 & 0.0320 & 0.0019 \\
  $27\le j<38$ & 4.816 & 0.3490 & 0.2541 & 0.1747 & 0.1271 & 0.0635 & 0.0317 & 0.0003 \\
  $38\le j<49$ & 4.822 & 0.3492 & 0.2541 & 0.1745 & 0.1271 & 0.0634 & 0.0317 & 0.0002 \\
  $49\le j<60$ & 4.822 & 0.3497 & 0.2538 & 0.1745 & 0.1268 & 0.0634 & 0.0318 & 0.0005 \\
  \midrule
  prediction & 4.819 & 0.3492 & 0.2540 & 0.1746 & 0.1270 & 0.0635 & 0.0317 & \\
  \bottomrule
 \end{tabular}
\end{table}

\subsection*{Consequences for the terminal branches}
Formally differentiating the octave law suggests the pointwise approximation
\begin{equation}\label{eq:pointwise}
 \frac{|\Pred(a)\cap[1,N]|}{N}\ \approx\ \kappa\,\frac{h(a)}{a},\qquad \kappa=\frac{2}{\log(4/3)},
\end{equation}
where $h$ stands for the density of the Syracuse law near $a$, approximated by $3^m\pi_m(a\bmod3^m)$. This cannot hold literally: by \Cref{rem:unbounded} the densities $3^m\pi_m$ do not converge uniformly, and for small $a$ the additive term in $3n+1$ is not negligible. Nevertheless it captures the main variation. For $2^{12}\le a<2^{15}$, $3\nmid a$, the correlation coefficient between the logarithms of $a\,|\Pred(a)\cap[1,N]|/(\kappa N)$ and of $3^m\pi_m(a\bmod 3^m)$ rises steadily from $0.69$ at $m=2$ to $0.945$ at $m=10$ (from $0.27$ to $0.71$ without logarithms), and the mean of the former quantity is $1.46$, against the value $\frac32$ predicted by the octave law (units make up two thirds of the integers and carry all the mass). The last two columns of \Cref{tab:delta} compare both sides of \eqref{eq:pointwise} for $a=q_k$. The residues of the $q_k$ modulo $9$ are periodic in $k$ with period $9$; those with $k\equiv2\pmod3$ run through $5,8,2$ (forms $\fe,\fa,\fc$, of Syracuse mass $\frac{4}{63},\frac{22}{63},\frac{16}{63}$) and those with $k\equiv1\pmod3$ through $1,4,7$ (forms $\fd,\fb,\ff$, of mass $\frac{8}{63},\frac{11}{63},\frac{2}{63}$). The exceptionally thin branch at $q_7=5461\equiv7\pmod9$ sits in the lightest of the six forms, and the heavy branches at $q_5$ and $q_8$ in the two heaviest.

The picture that emerges is as follows. The partition of the integers into terminal branches is governed by two factors: the factor $1/q_k$, which makes the first branch $q_2=5$ dominant, and a $3$-adic factor given by the Syracuse density at $q_k$, which modulates the remaining ones. The six forms thus appear twice: as the residue of the terminal power $2^{2k-1}$, which labels the classes, and as the residue of the last odd term $q_k$, which governs their sizes.

\begin{conjecture}[Terminal constants]\label{conj:alpha}
 For every $k\ge2$ with $3\nmid k$ the set $\Pred(q_k)$ has a natural density $\delta_k>0$, $\sum_k\delta_k=1$, and
 \[
  \lim_{N\to\infty}\Freq_{\fa}(N)=\alpha=\sum_{k\equiv2\ (3)}\delta_k,\qquad\lim_{N\to\infty}\Freq_{\fc}(N)=1-\alpha=\sum_{k\equiv1\ (3),\ k\ge4}\delta_k .
 \]
 Numerically $\delta_2\approx0.9377$, $\delta_4\approx0.02376$, $\delta_5\approx0.03796$ and $\alpha\approx0.9761$.
\end{conjecture}

The empirical frequencies still fluctuate in the fifth decimal place at $N=4\times10^9$ (\Cref{tab:freq}), so we do not claim more digits. We know of no closed form for $\alpha$; in view of \eqref{eq:predsplit}, the densities $w(a)$ of the sets $\Pred(a)$, if they exist, satisfy the infinite linear system $w(a)=w(2a)+\mathbf 1[a\equiv2\ (3)]\,w(\frac{2a-1}{3})$, of which \eqref{eq:pointwise} describes the expected behaviour at infinity.

\section{Open problems}\label{sec:open}

\begin{enumerate}[label=\textbf{Q\arabic*.}]
 \item \textbf{Cutoff.} Prove \Cref{conj:cutoff}, or at least improve either constant in $k\log_2\frac32\lesssim t_{\mathrm{mix}}^{(k)}\lesssim2k$.
 \item \textbf{Integrability.} Is $c_k=3^k\sum_x\pi_k(x)^2$ bounded, of linear growth, or of power growth? Equivalently, is the density of the Syracuse measure on $\Z_3$ square integrable?
 \item \textbf{The mode.} Is $-1$ the unique mode of $\pi_k$ for every $k$, and does $2^k\pi_k(-1)$ converge?
 \item \textbf{Octave law.} Prove \Cref{conj:octave} in any averaged form. For heights $y=N^\theta$ with $\theta>\log3/\log4$ and orbit points counted during the first $\log_2 N$ steps, a statement of this kind should follow from \Cref{thm:exact} and the Markov renewal theorem; the difficulty is to descend to fixed heights.
 \item \textbf{Terminal constants.} Do the densities $\delta_k$ exist? Is there a description of $\alpha$ other than as the solution of the infinite system above?
 \item \textbf{Other maps.} For the $qx+1$ maps ($q$ odd) the analogous chain acts on $\Z/q^k\Z$ with maps $r\mapsto r/2$, $r\mapsto(qr+1)/2$, and the argument proving \Cref{thm:spectrum} shows again that all non-trivial eigenvalues have modulus $\frac12$. How do the stationary laws and the terminal statistics compare?
\end{enumerate}

\section*{Computations and data availability}
All computations were carried out with short C and Python programs: an exhaustive first-passage sieve for $n\le4\times10^9$ (64-bit arithmetic with overflow check), exact iteration of $P_k$ in double precision for $k\le12$, and a Monte Carlo simulation with a fixed seed. The stationary laws for $k\le3$ and the characteristic polynomials for $k\le3$ were additionally verified in exact rational, respectively floating-point, arithmetic, and \Cref{thm:exact} was checked by enumeration for $k=2$, $j=10$. The programs accompanying this manuscript are available at \url{https://github.com/Sawon1234/Collatz_Automata}.




\end{document}